\documentclass[12pt, reqno]{amsart}
\usepackage{amsmath, amsthm, amscd, amsfonts, amssymb, graphicx, xcolor}
\usepackage[bookmarksnumbered, colorlinks, plainpages]{hyperref}
\usepackage{setspace}
\newtheorem{theorem}{Theorem}[section]

\newtheorem{proposition}[theorem]{Proposition}

\theoremstyle{definition}
\newtheorem{definition}[theorem]{Definition}

\theoremstyle{remark}

\numberwithin{equation}{section}

\begin{document}
\setcounter{page}{1}

\centerline{}

\centerline{}

\title[Negative Schwarzian derivative in Real One-Dimensional Dynamics]{On the Use of the Schwarzian derivative in Real One-Dimensional Dynamics}


\author[Bernardo San Martín, Felipe Correa]{Bernardo San Martín $^1$ and Felipe Correa$^2$$^{*}$}

\address{$^{1}$ Departamento de Matemáticas, Universidad Católica del Norte, Antofagasta, Chile.}
\email{\textcolor[rgb]{0.00,0.00,0.84}{sanmarti@ucn.cl}}

\address{$^{2}$ Departamento de Matemáticas, Universidad Católica del Norte, Antofagasta, Chile.}
\email{\textcolor[rgb]{0.00,0.00,0.84}{fcorrea02@ucn.cl}}




\begin{abstract}
In the study of properties within one-dimensional dynamics, the assumption of a negative Schwarzian derivative has been shown to be very useful. However, this condition may seem somewhat arbitrary, as it is not inherently a dynamical condition, except for the fact that it is preserved under iteration. In this brief work, we show that the negative Schwarzian derivative condition is not arbitrary in any sense but is instead strictly related to the fulfillment of the Minimum Principle for the derivative of the map and its iterates, which plays a key role in the proof of Singer’s Theorem.

\end{abstract} \maketitle

\section{Introduction}
The {\em Schwarzian derivative} appears in a wide range of mathematical topics, often in areas that seem unrelated at first glance \cite{MR2489717}. It was first formulated in 1869 by Hermann A. Schwarz in his work on conformal mappings. David Singer was the first to apply it to one-dimensional dynamics in 1978, using it to study $C^3$ maps from the unit interval to itself. An initial approach was made by D.J. Allwright in \cite{MR475982} who studied bifurcations of $C^3$ maps satisfying a certain property, denoted by $P$, which resembles the Schwarzian derivative. In later progress, it was found that maps with negative Schwarzian derivative also possess local properties that are useful for establishing certain distortion bounds, particularly when focusing on cross-ratios \cite{MR929092}.

At first glance, the negative Schwarzian derivative condition may appear somewhat arbitrary because it does not seem to be a dynamical condition. In this note, we show that this condition is not arbitrary in any sense, but rather strictly related to a sufficient condition that guarantees the fulfillment of the Minimum Principle for the derivative of the map and its iterates, which is the key point in the proof of Singer's Theorem. To the best of our knowledge, this is a simple and illustrative explanation of the use of the Schwarzian derivative in the context of one-dimensional differential dynamics.
\section{The Schwarzian Derivative}
A widely held view among mathematicians who have worked with the Schwarzian derivative is the sense of mystery surrounding its origin and the remarkable way in which it facilitates solving various problems in one-dimensional dynamics, driven only by the requirement to preserve the negative condition under iterations of the map. However, its precise connection to the dynamical properties of the map remains unclear. Let us briefly recall it. Consider an interval $I=(a,b)\subseteq \mathbb{R}$ and a $C^3$ map $f:I\to I$. If $f'(x)\neq 0$,
    the {\em Schwarzian derivative} of $f$ at $x$ is defined as
    $$Sf(x):=\frac{f'''(x)}{f'(x)}-\frac{3}{2}\left(\frac{f''(x)}{f'(x)}\right)^2.$$
We say that $f$ has a {\em negative Schwarzian derivative} on $I$ if $Sf(x)<0$ for all $x\in I$.

One of the main reasons the Schwarzian derivative is of interest in one-dimensional dynamics, as first observed by Singer, is its remarkable composition law, which follows directly from the chain rule
\begin{equation}\label{comp}
    S(h\circ g)(x)=Sh(g(x))\cdot(g'(x))^2+Sg(x).
\end{equation}
Consequently, if a map has a negative Schwarzian derivative, so do all its iterates. As we shall see, this property makes it a valuable tool in one-dimensional dynamics.
\section{The Minimum Principle}


\begin{definition}
    [The Minimum Principle in an interval]\label{MP}
     A map $g$ defined in an interval $J=[a,b]$ satisfies {\em the Minimum Principle on $J$} if for all $x\in (a,b)$
     $$|g(x)|>\min \{|g(a)|,|g(b)|\}.$$
\end{definition}

\begin{definition}[The Minimum Principle]
     A map $g$ defined on an interval $I$ satisfies the {\em The Minimum Principle} if it satisfies the Minimum Principle in all intervals $J\subset I$  where the map $g$ does not vanish.
\end{definition}

By definition, a map $g$ satisfies the Minimum Principle if and only if every local maximum is positive and every local minimum is negative.

In particular, for a differentiable map $f$ defined on an interval $I$, its derivative $f'$ satisfies the Minimum Principle if for any non-vanishing critical point $x\in I$ of $f'$, the quotient ${\displaystyle \frac{f'''(x)}{f'(x)}}$ is negative.

\section{Singer's Theorem}
An important result in one-dimensional differential dynamics was stated by Singer in \cite{MR494306}. This result shows that the negative Schwarzian derivative assumption imposes conditions on the type and quantity of periodic orbits that the map can possess.

Before stating the Theorem, we recall some necessary definitions. We say that $p$ is a \textit{periodic point} for a map $f$ if, for some integer $n$, $f^n(p)=p.$ Denote $O(p)=\{f^n(p);\,n\in\mathbb{Z}\}$ the {\em orbit} of $p$ under $f$. The {\em $\omega-$limit set} is the set of accumulation points of the sequence of forward iterates of a point in this orbit. The {\em basin} of a periodic point $p$ is the set of points whose $\omega-$limit set contains $p$. A periodic point $p$ is {\em attracting} if its basin contains an interval that contains $p$.
The {\em immediate basin} of a periodic point $p$ is the union of the connected components of its basin which contain a point from $O(p)$. Finally, we say that $c$ is a critical point of  $f$ 
 if $f'(c) = 0$. A critical point is called {\em non-degenerate} if $f''(c)\neq 0$.
\begin{theorem}[Singer's Theorem \cite{MR494306}]
    If $f:I\to I$ is a $C^3$ map with negative Schwarzian derivative, then the immediate basin of any attracting periodic point contains either a critical point of $f$ or a boundary point of $I$; each neutral periodic point is attracting; and there exists no interval of periodic points.
\end{theorem}
The key point in the proof of Singer's Theorem follows from the following proposition \cite{MR494306}.
\begin{proposition}
    If $Sf(x)<0$ for all $x\in I$,  then the function $f'$ cannot have either a positive local minimum value or a negative local maximum value.
\end{proposition}
In particular, the negative Schwarzian derivative condition for $f$, combined with the composition law in Eq. (\ref{comp}), implies that for all positive integers $n$, the derivative $(f^n)'$ satisfies the Minimum Principle. 

\section{The Schwarzian Derivative \& The Minimum Principle}


Our goal is to find a condition on $f$ (without using the Schwarzian derivative and its properties) such that the derivative $(f^n)'$ satisfies the Minimum Principle for all positive integers $n$. Unfortunately, deriving such a condition does not seem to be a straightforward task. Instead, we propose an alternative approach: finding a condition on $f$ such that, for a given positive integer $n$ and a non-vanishing critical point $x\in I$ of $(f^{n+1})'$, the quotient
$$\frac{(f^{n+1})'''(x)}{(f^{n+1})'(x)} $$
be negative. Indeed, as noted at the end of Section 3 this will guarantee the fulfillment of the Minimum Principle for $(f^{n+1})'$.

To achieve this, let $f$ be a differentiable map defined on an open interval $I$, and let $x\in I$ be a non-vanishing critical point of $(f^{n+1})'$; that is, a point $x\in I$ such that $(f^{n+1})'(x)\neq 0$ and $(f^{n+1})''(x)=0$.

By the chain rule, we have the following
\begin{equation}\label{der111}
            (f^{n+1})'(x)=(f^n)'(f(x))\cdot f'(x);
        \end{equation}
        \begin{equation}\label{der211}
            (f^{n+1})''(x)  =  (f^n)''(f(x))\cdot (f'(x))^2+(f^n)'(f(x))\cdot f''(x);
        \end{equation}
        and
        \begin{multline}\label{der311}
            (f^{n+1})'''(x)  =  (f^n)'''(f(x))\cdot (f'(x))^3+3\,(f^n)''(f(x))\cdot f'(x)\cdot f''(x)\\ +(f^n)'(f(x))\cdot f'''(x).
        \end{multline}
        Thus, from Eqs. (\ref{der111}) and (\ref{der311}), we obtain
        \begin{equation}\label{eq41}
            \frac{(f^{n+1})'''(x)}{(f^{n+1})'(x)} = \frac{(f^n)'''(f(x))}{(f^n)'(f(x))}\cdot (f'(x))^2 +3\,\frac{(f^n)''(f(x))\cdot f''(x)}{(f^n)'(f(x))} +\frac{f'''(x)}{f'(x)}.
        \end{equation}
      Since we have assumed that $(f^{n+1})''(x)=0$, it follows from Eq. (\ref{der211}) that
        \begin{equation}\label{eq511}
            (f^n)''(f(x))\cdot (f'(x))^2=-(f^n)'(f(x))\cdot f''(x).
        \end{equation}
        
        First, multiply Eq. (\ref{eq511})  by $f''(x)$  and rearrange terms, we obtain
        \begin{equation*}
            \frac{(f^n)''(f(x))\cdot f''(x)}{(f^n)'(f(x))}=-\left(\frac{f''(x)}{f'(x)}\right)^2.
        \end{equation*}
       Substituting this into the second term on the right-hand side of Eq. (\ref{eq41}) gives
        \begin{equation}\label{A1}
           \frac{(f^{n+1})'''(x)}{(f^{n+1})'(x)} =  \frac{(f^n)'''(f(x))}{(f^n)'(f(x))}\cdot (f'(x))^2 -3\,\left(\frac{f''(x)}{f'(x)}\right)^2+\frac{f'''(x)}{f'(x)}.
        \end{equation}

        Second, again from Eq. (\ref{eq511}), dividing by $((f^n)'(f(x)))^2$, we obtain
        \begin{equation*}
          \frac{(f^n)''(f(x))}{((f^n)'(f(x)))^2}\cdot (f'(x))^2=  -\frac{f''(x)}{(f^n)'(f(x))}.
        \end{equation*}
       Substituting this into the second term on the right-hand side of Eq. (\ref{eq41}) gives
\begin{multline}\label{B1}
    \frac{(f^{n+1})'''(x)}{(f^{n+1})'(x)}  =  \frac{(f^n)'''(f(x))}{(f^n)'(f(x))}\cdot (f'(x))^2-3\left(\frac{(f^n)''(f(x))}{(f^n)'(f(x))}\right)^2\cdot (f'(x))^2 +\frac{f'''(x)}{f'(x)}.
\end{multline}

Thus, by adding Eqs. (\ref{A1}) and (\ref{B1}), dividing by 2 and rearranging terms, we obtain
\begin{multline}\label{schwaff1}
    \frac{(f^{n+1})'''(x)}{(f^{n+1})'(x)}  = \left(\frac{(f^n)'''(f(x))}{(f^n)'(f(x))}-\frac{3}{2}\left(\frac{(f^n)''(f(x))}{(f^n)'(f(x))}\right)^2\right)\cdot (f'(x))^2\\ +\frac{f'''(x)}{f'(x)}-\frac{3}{2}\left(\frac{f''(x)}{f'(x)}\right)^2.
\end{multline}
Hence, assuming that $(f^{n+1})''(x)=0$ and that both expressions
\begin{equation}\label{schwa21}
    \frac{(f^n)'''(f(x))}{(f^n)'(f(x))}-\frac{3}{2}\left(\frac{(f^n)''(f(x))}{(f^n)'(f(x))}\right)^2
\end{equation}
and
\begin{equation}\label{schaw11}
\frac{f'''(x)}{f'(x)}-\frac{3}{2}\left(\frac{f''(x)}{f'(x)}\right)^2
\end{equation}
are negative, it follows from Eq. (\ref{schwaff1}) that the quotient $$\frac{(f^{n+1})'''(x)}{(f^{n+1})'(x)}$$ is negative, as desired.

Thus, for a given positive integer $n$, the Minimum Principle holds for $(f^{n+1})'$ provided that, for any non-vanishing critical point $x\in I$ of $(f^{n+1})'$, the expressions in Eqs. (\ref{schwa21}) and (\ref{schaw11}) are negative.

Note that the expressions in Eqs. (\ref{schwa21}) and (\ref{schaw11}) correspond to $S(f^n)(f(x))$ and $Sf(x)$, respectively. Therefore, for any non-vanishing critical point $x\in I$ of $(f^{n+1})'$, Eq. (\ref{schwaff1}) can be rewritten as
\begin{equation}\label{ber1}
    \frac{(f^{n+1})'''(x)}{(f^{n+1})'(x)}=S(f^n)(f(x))\cdot (f'(x))^2+Sf(x),
\end{equation}
which resembles the composition law given in Eq. (\ref{comp}), with $h=f^n$ and $g=f$. In fact, by the definition of the Schwarzian derivative for $f^{n+1}$ and Eq. (\ref{ber1}), it follows directly that, for any non-vanishing critical point $x\in I$ of $(f^{n+1})'$,
\begin{equation}\label{ber2}
    S(f^{n+1})(x)=S(f^n)(f(x))\cdot (f'(x))^2+Sf(x).
\end{equation}
On the other hand, a straightforward computation shows that Eq. (\ref{ber2}) holds for all $x\in I$, leading to a remarkable consequence: the negative Schwarzian derivative condition for $f$ is preserved by iterates. Therefore, the Minimum Principle for the iterates of $f$ is guaranteed only by the requirement that the expression in Eq. (\ref{schaw11}) is negative for all $x\in I$, which precisely corresponds to the negative Schwarzian derivative condition for $f$.
\section{Conclusion}
The negative Schwarzian derivative condition naturally emerges when looking for a condition on $f$ that guarantees that $(f^{n})'$ satisfies the Minimum Principle for all positive integers $n$. This principle has significant dynamical implications, particularly with respect to the type and number of positively oriented fixed points of $f$, as stated in Singer's Theorem. To the best of our knowledge, in the literature, there is no satisfactory and simple explanation about the use of the Schwarzian derivative in Real One-Dimensional Dynamics. Therefore, we believe that this brief note is valuable in order to shed light on this point.

\bibliographystyle{amsplain}
\bibliography{SchwDer.bib}

\end{document}